\documentclass[12pt]{amsart}

\usepackage{xy, graphicx, color, hyperref,array, mathtools}
\usepackage{hyperref}

\usepackage{amsmath}
\usepackage{amssymb}
\usepackage{amsfonts}
 \usepackage{bm}
\usepackage{mathrsfs}

\usepackage{tikz,verbatim}
 \usepackage{dsfont}

\makeatletter
\def\subsection{\@startsection{subsection}{3}%
  \z@{.9\linespacing\@plus.7\linespacing}{.1\linespacing}%
  {\normalfont\bfseries}}
\makeatother 
 
 \xyoption{all} 
\title[Divisibility of Degrees]{On the Divisibility of Degrees of Representations of  Lie Algebras}
 
 \author{Varun Shah, Steven Spallone }

 \newtheorem{thm}{Theorem}[section]
\newtheorem{c.intro}[thm]{Corollary}
\newtheorem{lemma}[thm]{Lemma}
\newtheorem{prop}[thm]{Proposition}
\newtheorem{cor}{Corollary}[thm]

\theoremstyle{definition}
\newtheorem{remark}[thm]{Remark}
 
\newtheorem{defn}[thm]{Definition} 
\newtheorem{example}[thm]{Example} 
\newcommand{\nc}{\newcommand}

\nc{\theorem}{\thm}

 \nc{\corollary}{\cor}
\nc{\mc}{\mathcal}
\nc{\mb}{\mathbb}
\nc{\mf}{\mathfrak}
\nc{\ul}{\underline}
\nc{\ol}{\overline}
\nc{\N}{\mb N}
\nc{\R}{\mb R}
\nc{\Z}{\mb Z}
\nc{\Q}{\mb Q}
\nc{\C}{\mb C}
\nc{\ms}{\mathscr}

\nc{\dmo}{\DeclareMathOperator}

\nc{\mat}[4]{
    \begin{pmatrix}
      #1 & #2 \\
      #3 & #4
    \end{pmatrix}
}

\dmo{\Ker}{Ker} \dmo{\val}{val} \dmo{\ord}{ord}
\dmo{\simp}{sc}
 \dmo{\wdf}{wdf}
 \dmo{\WDF}{WDF}

\dmo{\odd}{odd}
\dmo{\sgn}{sgn}
\nc{\beq}{\begin{equation*}}
\nc{\eeq}{\end{equation*}}
\nc{\half}{\frac{1}{2}}

\dmo{\Mod}{mod}

\dmo{\core}{core}
\dmo{\res}{res}
\dmo{\lin}{lin}
 \dmo{\vol}{vol}
\dmo{\Sp}{Sp}
\dmo{\SO}{SO}
\dmo{\PGL}{PGL}
\dmo{\SL}{SL}
 \dmo{\Spin}{Spin}
\dmo{\GSp}{GSp}
\nc{\la}{\lambda}
  \nc{\eps}{\varepsilon}
  
 \nc{\lip}{\langle}
 \nc{\rip}{\rangle}
\nc{\gm}{\gamma}

\dmo{\sd}{sd}
\dmo{\Res}{Res}
\dmo{\Ind}{Ind}
\dmo{\tr}{tr}
\dmo{\Sym}{Sym}
\dmo{\reg}{reg}
\dmo{\End}{End}
\dmo{\Hom}{Hom}
\dmo{\orth}{orth}
\dmo{\Prop}{Prop}

 \address{Indian Institute of Science Education and Research, Pune-411008,Maharashtra,India}
\email{shah.varun@students.iiserpune.ac.in}
\address{Indian Institute of Science Education and Research, Pune-411008,Maharashtra,India}
\email{sspallone@gmail.com}
\date{\today}

\keywords{Weyl dimension formula, natural density, integer-valued polynomials}
\subjclass{Primary 22E46, Secondary 05E10}
 
\begin{document}
\maketitle

 \begin{abstract} Let $\mf g$ be a reductive Lie algebra, and $m$ a positive integer. There is a natural density of  irreducible representations of $\mf g$, whose degrees are not divisible by $m$. For $\mf g=\mf{gl}_n$, this density decays exponentially to $0$ as $n \to \infty$. Similar results hold for simple Lie algebras and Lie groups, and there are versions for self-dual and orthogonal representations. 
 \end{abstract}

 \tableofcontents
 
 \section{Introduction}
  
Let $m$ be a positive integer. In his 1971 paper \cite{macdonald}, Macdonald gave a formula for the number of irreducible representations, or Specht modules, of the symmetric group $S_n$, whose degrees are not divisible by $m$. Dividing his formula by the number of Specht modules gives a density $0 \leq \mf d_m(S_n) \leq 1$ of Specht modules with this property. It is not hard to show (e.g., as in \cite{GPS}) that 
\begin{equation} \label{macd}
\lim_{n \to \infty} \mf d_m(S_n)=0.
\end{equation}
 In this paper we find analogous  asymptotic statements in the context of Lie groups and Lie algebras. 
 
 Let $\mf g$ be a reductive complex Lie algebra. Write $\Lambda$ for the weight lattice, and $\Lambda^+$ for the dominant weights. Then $\Lambda^+$ parametrizes the irreducible representations of $\mf g$. Write $\Lambda_m^+ \subseteq \Lambda^+$ for the subset of representations whose degrees are not divisible by $m$. 
What proportion of weights of $\Lambda^+$ lie in $\Lambda_m^+$?  For instance, when $\mf g=\mf{sl}_2$, there is exactly one irreducible representation of each degree, so the proportion should be $\frac{m-1}{m}$.
 
To address this question, we develop a notion of density for a subset $A$ of $\Lambda^+$, written $\mf d(A \mid \Lambda^+)$.  Our definition generalizes the concept of natural density for subsets of integers, described for instance in \cite[Chapter 11]{niven}. It should be interpreted as the probability that a randomly chosen member of
 $\Lambda^+$ lies in $A$. We put $\mf d_m(\mf g)=\mf d(\Lambda_m^+ \mid \Lambda^+)$; this is our answer to the above question.

 \begin{thm}\label{main.thm}We have
\beq
\frac{n!}{(mn)^n} \leq \mf d_m(\mf{gl}_n) \leq \omega(m)\exp \left(- \frac{n}{4m} \right).
\eeq
\end{thm}
   
(Recall that  $\omega(m)$ is the number of distinct prime factors of $m$.)
 
 We study these densities for classical Lie algebras, and obtain the following analogue of \eqref{macd}:

\begin{thm} \label{intro.simple.thm} Let $\mf g_j$ be a sequence of simple Lie algebras, whose dimensions tend to infinity as $j \to \infty$. Then for any positive integer $m$, the densities $\mf d_m(\mf g_j)$ decay exponentially to $0$ as $j \to \infty$.
\end{thm}

Similar asymptotics hold when we restrict our attention to Lie group representations, self-dual representations, and orthogonal representations.
 
 Let us give some idea of the proof. Since the function  $\Lambda \to \Z$ defined by
 \beq
\la \mapsto \deg \pi_{\la}
 \eeq
 is an integer-valued polynomial, its values mod $m$ are periodic. (See Lemma \ref{manyvariate}.)
  This allows us to find the density of $\Lambda_m^+$ in $\Lambda$ by counting its points in a suitable fundamental domain.
 
 The notion of density is introduced and developed in Section \ref{periodic.sets}. In Section \ref{integer.poly.section} we review the theory of integer-valued polynomials, and study their periodicity mod $m$. 
  
 Character degrees of general linear and classical Lie algebras are divisible by 
 certain Vandermonde determinants. In Section \ref{vandy}, we show that values of such polynomials tend to be highly divisible.   
 Theorem \ref{main.thm} essentially follows from this.
 
 We apply our results to irreducible representations of Lie algebras and Lie groups in Section \ref{rep.app.section}, and deduce Theorem \ref{intro.simple.thm}. 
 We sketch similar results  for self-dual and orthogonal representations in Section \ref{SDR.sec}.     
   \bigskip
   
\textbf{Acknowledgements.} 
The authors would like to thank Rohit Joshi for helpful discussions.

\section{Density of Periodic Sets}\label{periodic.sets}

\subsection{Vectorial Natural Density}\label{def.and.prop}
Let $L$ be a free $\Z$-module of finite rank, and $A$ a subset of $L$. Put $V=L \otimes_\Z \R$.  Given a norm $N$ on $V$, write
\begin{equation} \label{density.N}
\mf d_N(A \mid L)=\lim_{r \to \infty}\frac{ \# \{ a \in A : N(a)<r\}}{\# \{ v \in L: N(v)<r\}},
\end{equation}
provided this limit exists. 

We say a positive integer $\varpi$ is a \emph{period} of $A$ when $A+\varpi L = A$. When $A$ has a period, we say that $A$ is \emph{periodic}. Suppose $\varpi$ is a period of $A$, and $\mc B=\{v_1, \ldots, v_n\}$ is a basis of $L$. Put
\beq
\mc F_{\mc B}=\left\{ \sum_i c_iv_i  \mid c_i \in \Z \text{ and } 0 \leq c_i<\varpi \right \};
\eeq
it is a transversal for $L/\varpi L$.

The following theorem shows that density is a sensible notion for periodic subsets, and gives a expression to compute it.

\begin{thm} \label{mega.density.thm} Suppose that $A$ is a periodic subset of $L$. Then for any choices of norm $N$, period $\varpi$, and basis $\mc B$, we have
\beq
\mf d_N(A \mid L)=\frac{|\mc F_{\mc B} \cap A|}{|\mc F_{\mc B}|}.
\eeq
In particular it is independent of these choices. We hence define
\beq
 \mf d(A \mid L)=\mf d_N(A \mid L).
 \eeq
When $A$ is a sublattice of $L$, we have
\beq
\mf d(A \mid L)=\frac{1}{[L:A]}.
\eeq
If $x \in L$, then
\beq
\mf d(A +x \mid L)=\mf d(A \mid L).
\eeq
\end{thm}

This will be proved by the end of this section. 

\begin{remark} See \cite[Section 2]{maze}, and the citations therein, for earlier formulations of  natural density for $\Z^n$.
\end{remark}

\begin{remark} The limit in \eqref{density.N} may not exist in general, even when $n=1$. See \cite[Exercise 9, p.~ 475]{niven}.
\end{remark}

\subsection{Proof of Formula}

In this subsection we prove the first statement of Theorem \ref{mega.density.thm}. 

\begin{lemma} \label{shift.fund.lemma} 
For any $x \in  L$, we have
\beq
|\mc F_{\mc B} \cap A|= |(\mc F_{\mc B}+x) \cap A|.
\eeq
\end{lemma}

\begin{proof} 
Reducing $\mc B$-coordinates modulo $\varpi$ maps $L$ to $\mc F_{\mc B}$; we write this as $v \mapsto v \% \varpi$.
It is elementary to see that the restriction of this map to $A$ gives a bijection
\beq (\mc F_{\mc B}+x)\cap A \overset{\sim}{\to} \mc F_{\mc B}\cap A. \eeq    \end{proof}

\begin{lemma} \label{convex.lemma} 
        Let $K$ be a compact, convex subset of  $V$ with nonempty interior. 
       Then 
        \beq
            \lim\limits_{r \to \infty} \frac{|rK \cap L|}{r^n} = M_V
        \eeq
        for some constant $M_V > 0$.
        
\end{lemma}
    
\begin{proof}     
   This is a basis-free version of   \cite[p.~120, Equation (2) and Theorem 4]{gruber}.
   (The constant $M_V$ is the Euclidean volume of $K$ with respect to a basis of $L$.)
    \end{proof}

\begin{prop} \label{fund.domain.prop}
For any period $\varpi$ of $A$, any norm $N$ of $V$, and any basis $\mc B$ of $L$, we have
\beq
\mf d_N(A \mid L)=\frac{| \mc F_{\mc B} \cap A|}{|\mc F_{\mc B}|}.
\eeq
\end{prop}
 
\begin{proof} 
            Translates of $\mc F_{\mc B}$ by vectors in $\varpi L$ partition $L$. We call each such translate a \textit{tile},  and
            say a subset of $L$ is \emph{tiled}, when it is a disjoint union of tiles.    
             Let 
            \beq
            C_r = \{x \in \varpi L \mid (x + \mc F_{\mc B}) \cap B_N(r) \neq \varnothing\},
            \eeq
          so that $S_r = \bigcup\limits_{x \in C_r} (x + \mc F_{\mc B})$ is the smallest tiled set containing $B_N(r) \cap L$.  If $d$ is the $N$-diameter of $\mc F_{\mc B}$ then $S_r \subseteq B_N(r+d)$. Indeed, let $v$ be an element of $S_r$, then it is in some tile, and this tile contains a $u \in B_N(r)$. Whence 
            \beq N(v) \leq N(u) + N(v-u) \leq r + d.\eeq
 In particular, we have $B_N(r) \cap L \subseteq S_r \subseteq B_N(r+d) \cap L$, which implies
            \beq
            |B_N(r) \cap L| \leq |S_r| \leq |B_N(r+d) \cap L|.\eeq
 Dividing by $|B_N(r) \cap L|$   gives  
            \beq
            \begin{split}
            1 &\leq \lim\limits_{r \to \infty} \frac{|S_r|}{|B_N(r) \cap L|} \\ 
            &\leq \lim\limits_{r \to \infty} \frac{|B_N(r+d) \cap L|}{|B_N(r) \cap L|} \\
             &= \lim_{r \to \infty} \frac{(r+d)^n}{r^n} \text{ (by Lemma \ref{convex.lemma})} \\           
            &= 1. \\
            \end{split}
\eeq
            Therefore
            \beq 
            \lim\limits_{r \to \infty} \frac{|S_r|}{|B_N(r) \cap L|} = 1.
            \eeq
 As before we have $S_{r-d} \subseteq B_N(r) \cap L \subseteq S_r$,  and the inclusions are preserved when we consider the intersections of these sets with $A$. Thus 
               \beq
               |S_{r-d} \cap A| \leq |B_N(r) \cap A| \leq |S_r \cap A|.
               \eeq
 Dividing everything by $|B_N(r) \cap L|$ and taking   $r \to \infty$ gives
            
            \beq \lim\limits_{r\to\infty} \frac{|S_{r-d} \cap A|}{|S_{r-d}|} \leq \lim\limits_{r\to\infty} \frac{|B_N(r) \cap A|}{|B_N(r) \cap L|} \leq \lim\limits_{r\to\infty} \frac{|S_{r} \cap A|}{|S_{r}|}.\eeq
            
            Since $S_r$ is the disjoint union of $|C_r|$ many tiles, $|S_r| = |C_r|\cdot |\mc F_{\mc B}|$. By Lemma \ref{shift.fund.lemma} each tile has the same number of points in $A$ so $|S_r \cap A| = |C_r| \cdot |\mc F_{\mc B} \cap A|$. Hence for any $r > 0$, 
            \beq
            \frac{|S_r \cap A|}{|S_r|} = \frac{|\mc F_{\mc B} \cap A|}{| \mc F_{\mc B}|}.
            \eeq         
          By the Squeeze Theorem we deduce 
           \beq
            \mathfrak d_N(A \mid L)= \frac{|\mc F_{\mc B} \cap A|}{| \mc F_{\mc B}|}.
            \eeq

        \end{proof}
 
\begin{cor} \label{interesting.problem} When $A$ is periodic, then $\mf d_N(A \mid L)$ exists and is   independent of $N$. The fraction
\beq
\frac{|\mc F_{\mc B} \cap A|}{|\mc F_{\mc B}|}
\eeq
is independent of $\varpi$ and $\mc B$.
\end{cor}

 To prove the last assertion in Theorem \ref{mega.density.thm}, note that 
 \beq |(A + x) \cap \mc F_{\mc B}| = |A \cap (\mc F_{\mc B} - x)| \eeq
 because the two sets differ by a translation of $x$.  Now apply Proposition \ref{fund.domain.prop} along with Lemma \ref{shift.fund.lemma}.
 
 \subsection{Sublattices}
In this section, we check that the density of a sublattice is equal to the reciprocal of its index. We also relate the density of a periodic subset of a sublattice $L'$  to its density in $L$.

Let $L$ be a free $\Z$-module of finite rank, and $L'$ a sublattice of full rank. By the theory of finitely generated abelian groups (\cite[Chapter 12, Theorem 4]{dummit}), there is a basis $v_1, \dots, v_n$ of $L'$ and positive integers $a_1, \dots, a_n$ so that $a_1v_1, \dots, a_nv_n$ is a basis of $L$ and              \beq
            a_1 \mid a_2 \mid \cdots \mid a_n.
            \eeq
        \begin{lemma}
          The sublattice $L'$ is an $a_n$-periodic subset of $L$ with density 
            \beq
            \mf d(L' \mid L) = \frac{1}{[L: L']}.
            \eeq
        \end{lemma}

        \begin{proof}
            
          Since
             $a_nL \subseteq L'$, it is plainly an $a_n$-periodic subset of $L$.  

            An element $v = \sum_i c_iv_i$, with $0 \leq c_i < a_n$, lies in $L'$ whenever each $a_i \mid c_i$. So there are $(a_n)^n/(a_1\cdots a_n)$ such points. This means that 
            \beq
            \mf d(L' \mid L) = \frac{1}{a_1\cdots a_n} = \frac{1}{[L:L']}.
            \eeq
        \end{proof}

        Let $V = L' \otimes_\Z \R$. We may identify $L$ with the $\Z$-span of the vectors $a_iv_i \otimes a_i^{-1} \in V$. 

        \begin{prop} \label{sublatticey}
            Let $L', L$ be as above and $A$ a periodic subset of $L$. Then 
            \beq
            \mf d(A\cap L' \mid L') \leq [L:L']\mf d(A \mid L).
            \eeq
        \end{prop}

        \begin{proof} 
            Let $N$ be a norm on $V$. Then for $A' = A\cap L'$, we have
            \beq
                \begin{split}
                    \mf d(A' \mid L') &= \lim\limits_{r \to \infty} \frac{|B_N(r) \cap A'|}{|B_N(r) \cap L'|} \\
                    &\leq \lim\limits_{r \to \infty} \frac{|B_N(r) \cap A|}{|B_N(r) \cap L'|} \\
                    &= \lim\limits_{r \to \infty} \frac{|B_N(r) \cap A|}{|B_N(r) \cap L|} \cdot  \lim\limits_{r \to \infty}\frac{|B_N(r) \cap L|}{|B_N(r) \cap L'|} \\
                    &= \mf d(A \mid L) [L : L'],\\
                \end{split}
            \eeq

            as required.
        \end{proof}

\subsection{Cones}
For our application to representation theory, we will need a ``cone'' version of density.
Again let $L$ be a lattice, and fix a basis $\mc B$ of $L$.
 The set $C$ of nonnegative linear combinations of $\mc B$ is a cone in $V$. 
Write $L^+=L \cap C$ and $A^+=A \cap L^+$.

\begin{defn} Given a norm $N$ on $V$, put
\beq
\mf d_N(A^+ \mid L^+)=\lim_{r \to \infty}\frac{ \# \{ a \in A^+ : N(a)<r\}}{\# \{ v \in L^+: N(v)<r\}}.
\eeq
\end{defn}

\begin{prop} \label{cone.version} When $A$ is periodic, the density $\mf d_N(A^+ \mid L^+)$ is independent of $N$, and
\beq
 \mf d(A^+ \mid L^+):=\mf d_N(A^+ \mid L^+)=\mf d(A \mid L).
\eeq
\end{prop}

\begin{proof} This follows as in the proof of Proposition \ref{fund.domain.prop}, replacing $L$ and $A$ everywhere with $L^+$ and $A^+$. (Note that $\mc F_{\mc B} \subset L^+$.)
\end{proof}

\subsection{A Nonperiodic Counterexample} 
In this section, we give an example to clarify that the hypothesis of periodicity in Theorem \ref{mega.density.thm} cannot be removed.   Let $L = \Z^2$ and take
    \beq A = \{(x, y) \in L \mid xy \geq 0\},\eeq
  namely the set of integer vectors in the first and third quadrant of $\R^2$. It is not periodic.  Consider two norms
    \beq N_0((x, y)) = \max(|x|, |y|) \eeq
    and
    \beq N_1((x, y)) = \max (|x|, |x+y|) \eeq
    on $\R^2$.

    From Lemma \ref{convex.lemma}, we see that for both $N = N_0$ and $N=N_1$, the quantity $ \mf d_N(A \mid L)$ is the area of the intersection of the first and third quadrants with $B_N(1)$.       Therefore
        \beq
            \frac{1}{2} = \mf d_{N_0}(A \mid L) \neq \mf d_{N_1}(A \mid L) = \frac{1}{4}.
        \eeq

\begin{remark} It seems an interesting question, to characterize the subsets of lattices which have a well-defined density (meaning independent of the norm), as in Corollary \ref{interesting.problem}.
        \end{remark}

\section{Integer-valued Polynomials} \label{integer.poly.section}

Congruence classes of binomial coefficients satisfy a certain periodicity. More precisely, let $p$ be a prime, and let $r,s,k$ be positive integers with $p^s \leq k <p^{s+1}$.
According to \cite[Theorem 4.8]{fray}, for all integers $n$ we have
\begin{equation} \label{FRA.thm}
\binom{n+p^{s+r}}{k} \equiv \binom{n}{k} \mod p^r.
\end{equation}

Let $m$ be a positive integer, and $f: \Z^n \to \Z$. 
By considering values of $f$ modulo $m$ we obtain a function 
\beq 
f \mod m: \Z^n \to \Z/m\Z.
\eeq
When there is a positive integer $\varpi$ so that 
\beq
f(v+\varpi w) \equiv f(v) \mod m
\eeq
for all $v,w \in \Z^n$, we say that $\varpi$ is an \emph{$m$-period} of $f$. 
 
 Write $\mb Q_{\mb Z}[x_1, \ldots, x_n]$ for the subring in $\mb Q [x_1, \ldots, x_n]$ of polynomials which take integer values on integers. For example $\half x(x+1) \in \mb Q_{\mb Z}[x]$. Given  $f \in \Q[x_1, \dots, x_n]$, denote by $\deg_{x_i}f$ the degree of $f$ in the variable $x_i$ and define $\deg_{\bullet}(f)= \max_{i}( \deg_{x_i} f)$. For instance,
 \beq
  \deg_\bullet(x^2y^2+xy)=2.
  \eeq
 
 The $\Z$-algebra $\mb Q_{\mb Z}[x]$ is freely generated by the binomial coefficients $\binom{x}{d}$ for $d=1,2,3, \ldots$. (See, for example, \cite{chabert}.) The following refinement
 can be proved by similar means (induction on the degree): 
 \begin{lemma} \label{refinement} A polynomial $f \in \Q_{\Z}[x_1, \dots, x_n]$ is in the subring of $\Q_{\Z}[x_1,\ldots, x_n]$ generated by $\dbinom{x_i}{k}$, with $k \leq \deg_\bullet(f)$.
\end{lemma}

 For a prime power $q = p^s$, put
\beq
r_q(f) = p^{\left\lfloor\log_p \deg_{\bullet}(f)\right\rfloor + s}.
\eeq

\begin{prop} \label{r.formula}The function $f \mod q$ is periodic with period $r_q(f)$.
\end{prop}

\begin{proof} 
This follows from \eqref{FRA.thm} and Lemma \ref{refinement}.
\end{proof}

Note the evident  bounds
\begin{equation} \label{evident.r}
\frac{q}{p}\deg_{\bullet}(f) \leq r_q(f) \leq q \deg_{\bullet}(f).
\end{equation}

 When $m$ is not a prime power, we can use the following straightforward observation:
 
 \begin{lemma} \label{manyvariate} Let $f \in \Q_\Z[x_1, \dots, x_n]$, and $m_1,m_2$ relatively prime integers.  If $\varpi_i$ is an $m_i$-period of $f$ for $i=1,2$, then $\varpi_1 \cdot \varpi_2$ is an $m_1m_2$-period of $f$.
 \end{lemma}

\section{Vandermonde Determinants}\label{vandy}
Let $L$ be a free $\Z$-module of finite rank $n$, and $\mc B=\{e_1,\ldots, e_n\}$ a basis of $L$. Write $\{e_1^*,\ldots, e_n^*\}$ for the dual basis of $L^\vee=\Hom(L,\Z)$. 
For $R=\Z$ or $\Q$, write $R[L^\vee]$ for the symmetric $R$-algebra $\Sym_{R}(L^\vee)$; one thinks of $R[L^\vee]$ as the ring of integer or rational polynomials on $L$, respectively. We associate to $\mc B$ the Vandermonde polynomial 
\beq
\mc V_{\mc B} =\prod_{i<j} (e_i^*-e_j^*) \in \Z[L^\vee].
\eeq
Write $(\mc V_{\mc B})$ for the ideal in $\Q[L^\vee]$ generated by $\mc V_{\mc B}$.   
 
	\begin{lemma} \label{main.bound.lemma} 
  Let $m$ be a positive integer, and $f \in (\mc V_{\mc B}) \cap \Q_\Z[L^\vee]$. Let $\varpi$ be an $m$-period of $f$, and
  put
  \beq 
  \mc F=\left\{\sum_{i=1}^n c_i e_i \mid 0 \leq c_i < \varpi \right\}.
  \eeq
 Then for $n \geq 2$ we have
            \beq \frac{\#\{\mathbf{x} \in \mc F : m \nmid f(\mathbf{x}) \}}{|\mc F|} \leq \exp{\left(-\frac{n^2}{4\varpi}\right)}.
            \eeq
	\end{lemma}	  	

        \begin{proof}
            The polynomial $\mc V_{\mc B}$ is skew-symmetric, in the sense that $\mc V_{\mc B}(\mathbf{x})$ vanishes whenever the $\mc B$-coordinates of $\mathbf{x}$ are not distinct. Since $f$ is a multiple of $\mc V_{\mc B}$, it is also skew-symmetric.              
       We have
          \beq
            \begin{split}
            \frac{\#\{\mathbf{x} \in \mc F : m \nmid f(\mathbf{x}) \}}{|\mc F|}    & \leq \frac{\#\{\mathbf{x} \in \mc F \mid \text{all coordinates of } \mathbf{x} \text{ are distinct}\}}{\varpi^n} \\
             &=    \frac{(\varpi)_n}{\varpi^n}=  \frac{\varpi}{\varpi}\frac{\varpi-1}{\varpi}\cdots \frac{\varpi-n+1}{\varpi}, \\
             \end{split}
             \eeq
           where $(\varpi)_n$ is the \emph{falling factorial} $\frac{\varpi !}{(\varpi - n)!}$.  
          Suppose first $\varpi \geq n$, so that each term in the product is positive.  By the  inequality of arithmetic and geometric means, we deduce
            \begin{align*}
                \frac{(\varpi)_n}{\varpi^n} &\leq \left(\frac{\varpi + (\varpi-1) + \cdots + (\varpi - n+1)}{n\varpi}\right)^n \\
                &= \left(1 - \frac{n-1}{2\varpi}\right)^n \\
                &\leq \left(1 - \frac{n}{4\varpi}\right)^n,
            \end{align*}
           the last inequality valid for $n \geq 2$.   For $x \leq 1$, we have the estimate
            \beq 0 \leq 1-x \leq e^{-x}.\eeq
Thus, as $4\varpi> n$, we have
            \beq\frac{(\varpi)_n}{\varpi^n} \leq \exp{\left(-\frac{n^2}{4\varpi}\right)}\eeq
            for all $n \geq 2$.

            On the other hand, if $\varpi < n$, then $(\varpi)_n/\varpi^n = 0$. So,   for all $n \geq 2$, we have
            \beq
            \frac{\#\{\mathbf{x} \in \mc F  : m \nmid f(\mathbf{x}) \}}{|\mc F|} \leq \frac{(\varpi)_n}{\varpi^n} \leq \exp{\left(-\frac{n^2}{4\varpi}\right)},
            \eeq
            as required.
        \end{proof}

        Suppose $m$ has prime factorization $m = p_1^{\alpha_1}\cdots p_r^{\alpha_r}$, so $m = q_1\cdots q_r$ where $q_i = p_i^{\alpha_i}$.

        \begin{prop}\label{cmi.fan}
            Let   $f \in (\mc V_{\mc B}) \cap \Q_\Z[L^\vee]$. For each $i$ pick a $q_i$-period $\varpi_{q_i}$ of $f$. Put
                \beq
 \bm{\varpi}   =\max_i (\varpi_{q_i})
                \eeq
            and let
                \beq
                    U_m(f) = \{\mathbf{x} \in L  : m \nmid f(\mathbf{x}) \}.
                \eeq
            Then for $n \geq 2$ we have
                \beq 
                    \mf d(U_m(f) \mid L) \leq \omega(m)\exp{\left(-\frac{n^2}{4\bm{\varpi} }\right)}.
                \eeq
        \end{prop}

        \begin{proof}
         
         Let $\varpi_m=\prod_i \varpi_{q_i}$; this is an $m$-period for $f$ by Lemma \ref{manyvariate}. So $U_m(f) \subset L$ has period $\varpi_m$. 
         
                      Given an $\mathbf{x} \in L$, clearly $m \nmid f(\mathbf{x})$ if and only if $q_i \nmid f(\mathbf{x})$ for some $i$.   By the union bound, we have:

            \begin{align*}
                \frac{\#\{\mathbf{x} \in \mc F  : m \nmid f(\mathbf{x}) \}}{|\mc F|} &\leq \sum\limits_{i} \frac{\#\{\mathbf{x} \in \mc F : q_i \nmid f(\mathbf{x})  \}}{|\mc F|} \\
                &\leq \sum_{i} \exp\left(-\frac{n^2}{4\varpi_{q_i}}\right) \text{ (by Lemma \ref{main.bound.lemma})}\\
                &\leq \omega(m) \exp{\left(-\frac{n^2}{4 \bm{\varpi}  }\right)}.
            \end{align*}

        The statement then follows from Proposition \ref{fund.domain.prop} with $A=U_m(f)$.
        \end{proof}
 

        \begin{cor}
            Let $m$ be as above, $f \in (\mc V_{\mc B}) \cap \Q_\Z[L^\vee]$ and $L'$ a sublattice of $L$. For each $i$, let  $\varpi_{q_i}$ be a $q_i$-period of $f$. Put
                \beq
                    \bm{\varpi} =\max_i (\varpi_{q_i})
                \eeq
            and let
                \beq
                    U'_m(f) = \{\mathbf{x} \in L'  : m \nmid f(\mathbf{x}) \}.
                \eeq
            Then for $n \geq 2$ we have
                \beq 
                    \mf d(U'_m(f) \mid L) \leq \omega(m)[L: L']\exp{\left(-\frac{n^2}{4\bm{\varpi} }\right)}.
                \eeq
        \end{cor}

        \begin{proof}
            Apply Propositions  \ref{sublatticey} and \ref{cmi.fan}.
        \end{proof}

\section{Application to Representations} \label{rep.app.section}

Now we apply the previous work to representations. A good general reference for this section is \cite{goodman}.
Every Lie group and Lie algebra considered in the following sections is complex.

\subsection{Density}

Let $\mf g$ be a reductive Lie algebra. 
Picking a Cartan subalgebra $\mf t$  of $\mf g$, one has a weight lattice  $\Lambda$. By picking a suitable ordering, one has a subset $\Lambda^+$ of dominant weights, which parametrizes the irreducible representations of $\mf g$. Write $\pi_\la$ for $\la \in \Lambda^+$. 

The Weyl dimension formula gives a polynomial function $\mathscr D=\mathscr D_{\mf g}: \Lambda \to \Q$, which restricts to
\beq
\mathscr D(\la)=\deg \pi_\la
\eeq
for $\la \in \Lambda^+$. By \cite[Proposition 25]{joshi}, $\ms D$ takes values in $\Z$, so 
$\mathscr D \in \Q_{\Z}[\Lambda^\vee]$. 
Let us set
\beq
\Lambda_m= \{ \la \in \Lambda : m \nmid \mathscr D(\la) \},
\eeq
so that
 \beq
  \Lambda_m^+ = \Lambda_m \cap \Lambda^+=\{ \la \in \Lambda^+ : m \nmid  \deg \pi_\la  \}.
  \eeq
 
\begin{defn} Let
\beq
\mf d_m(\mf g)=\mf d( \Lambda_m^+ \mid \Lambda^+).
\eeq 
\end{defn}
We interpret this as the proportion of irreducible representations of $\mf g$ whose degrees are not multiples of $m$.
By Proposition \ref{cone.version}, we have $\mf d_m(\mf g)=\mf d( \Lambda_m \mid \Lambda)$.

\subsection{The Lie Algebra $\mf{gl}_n$} \label{gln}

For $\mf g=\mf{gl}_n$ we may identify $\Lambda$ with $\Z^n$, and $\Lambda^+$ with the set of those 
\beq
\la=(\la_1, \ldots, \la_n) \in \Lambda
\eeq
with 
\beq
\la_1 \geq \cdots \geq \la_n \geq 0.
\eeq
(See \cite[Section 5.5.4]{goodman}.)  We put
\beq
\mathscr D_{n}(\la_1, \ldots, \la_n)= \mathscr D_{\mf {gl}_n}(\la)= \prod_{1 \leq i<j \leq r} \frac{\la_i-\la_j+j-i}{j-i}.
\eeq
Now put $\rho=(n-1,n-2, \ldots, 0)$ and consider
\beq
\mf f_n(\la)=\ms D_n(\la-\rho)= \prod_{1 \leq i<j \leq r} \frac{\la_i-\la_j}{j-i},
\eeq
so that $\mf f_n(x_1,\ldots, x_n) \in  (\mc V_n) \cap \Q_\Z[x_1, \dots,x_n]$.
Since $\deg_{x_i} \mf f_n=n-1$, we have $\deg_\bullet \mf f_n=n-1$. Let $\varpi_m$ be an $m$-period of $\mf f_n$. Put
\beq
\mc F_n =\{0, \ldots, \varpi_m-1\}^n.
\eeq

\begin{prop}  We have  
\beq
  \mf d_m( \mf{gl}_n)= \frac{ \# \{ {\bf x} \in \mc F_n : m \nmid \mf f_n({\bf x})   \}}{(\varpi_m)^n}.
\eeq
\end{prop}

\begin{proof} This is because
\beq
\begin{split}
\mf d_m( \mf{gl}_n) &= \mf d(\Lambda_m^+ \mid \Lambda^+) \\
				&= \mf d(\Lambda_m \mid \Lambda) \\
				&=\mf d(\{ {\bf x} \in \Lambda : m \nmid \ms D_n({\bf x})\} \mid \Lambda) \\
				&=\mf d(\{ {\bf x} \in \Lambda : m \nmid \ms D_n({\bf x}-\rho)\}-\rho \mid  \Lambda) \\
				&=\mf d(\{{\bf x} \in \Lambda : m \nmid \mf f_n({\bf x})\} \mid \Lambda), \\
				\end{split}
				\eeq
				and this may be computed using $\mc F_n$ by Theorem \ref{mega.density.thm}.
				\end{proof}

\begin{example} 
For $\mf g=\mf {gl}_2$, we have $\mf f_2(x_1,x_2)= x_1-x_2$. This gives $\varpi_m(\mf f_2)=m$ and $\mc F_2=\{0, \ldots, m-1\}^2$, so that
\beq
\mf d_m(\mf g)= \frac{m-1}{m}.
\eeq
For $n=3$ we have
\beq
 \mf f_3(x_1,x_2,x_3)=\frac{(x_1-x_2)(x_1-x_3)(x_2-x_3)}{6}.
\eeq
Let us take $m=2$, so that $\varpi_2(\mf f_3)=4$, and $\mc F_3=\{0,1,2, 3\}^3$. A small computation with $\mc F_3$ gives  
\beq
\mf d_2(\mf g)=\frac{3}{8}.
\eeq
\end{example}
 
We can continue this way. Here is   a table of values as $n$ grows:
\bigskip
\begin{center}
 \begin{tabular}{c|c| c|c|c|c|c|c|c| }
   & $1$ & $2$ & $3$ & $4$ & $5$  & $6$  & $7$ & $8$  \\
 \hline
 &&&&&&&& \\
  $m=2$ & $1$ & $ \frac{1}{2}$ & $\frac{3}{8}$ & $\frac{3}{32}$ & $\frac{15}{128}$ & $\frac{45}{1024}$ &  $\frac{315}{16384}$ & $\frac{315}{131072}$ \\
 &&&&&&&&  \\
  $m=3$ & $1$ & $\frac{2}{3}$ & $\frac{2}{9}$ & $\frac{8}{27}$ & $\frac{40}{243}$ & $\frac{80}{2187}$ & $\frac{560}{19683}$ 
  & $\frac{4480}{531441}$\\
  &&&&&&&&   \\
  \hline
 \end{tabular}
 
 \end{center}

\subsection{Proof of Theorem \ref{main.thm}}\label{gl.thm.proof}

 To establish the upper bound,
we apply Propositions \ref{cone.version}  and \ref{cmi.fan}, with $L=\Lambda$ and $L^+=\Lambda^+$, using the periods $r_q(\mf f_n)$ from Proposition \ref{r.formula}. 
Thus 
\beq
\mf d_m(\mf{gl}_n) \leq \omega(m)\exp\left(-\frac{n^2}{4 \bm{\varpi} } \right).
\eeq
 For each $i$ we have
 \beq
 r_{q_i}(\mf f_n) \leq \deg_\bullet(\mf f_n) \cdot q_i
 \eeq
 by \eqref{evident.r}. Since $\deg_\bullet(\mf f_n)=n-1$, we may simply write
 \beq
 \bm{\varpi}  =\max_i r_{q_i}(\mf f_n) \leq mn,
 \eeq
 and so 
 \beq
\mf d_m(\mf{gl}_n) \leq \omega(m) \exp\left(- \frac{n}{4m} \right).
\eeq
 This gives the upper bound of Theorem \ref{main.thm}.

For the lower bound, let $p$ be a prime divisor of $m$. From Proposition \ref{r.formula}, $\varpi = r_p(\mf f_n)$ is a $p$-period for $\mf f_n$. Because $\deg_{\bullet}(\mf f_n) = n-1$, we have
    \beq n-1 \leq \varpi \leq (n-1)p \leq mn \eeq
    from  \eqref{evident.r}.

    There are $n!$ permutations of the coordinates of $\rho$, written $\sigma \cdot \rho$ for $\sigma \in S_n$. For each, we have
    \beq
    \mf f_n(\sigma \cdot \rho) = \sgn(\sigma)= \pm 1,
    \eeq
      and therefore
    \beq
        \mf d_p(\mf{gl}_n) \geq \frac{n!}{\varpi^n} \geq \frac{n!}{(mn)^n}.
    \eeq
Since a multiple of $m$ is also a multiple of $p$, we deduce
    \beq \mf d_m(\mf{gl}_n) \geq \mf d_p(\mf{gl}_n) \geq \frac{n!}{(mn)^n}. \eeq
This proves the desired lower bound. $\square$

\subsection{The Lie Algebra $\mf{sl}_n$ } \label{slla.subsec}

In this section we treat the simple Lie algebra $\mf g'=\mf{sl}_n$, by relating it to $\mf g=\mf{gl}_n$. 
We drop the subscript `$n$', and  use prime notation (e.g., $\Lambda',\ms D'$) for the objects related to $\mf g'$, and maintain the earlier unprimed notation for objects related to $\mf g$.

   \begin{prop} \label{special.prop} For a positive integer $m$ we have
 \beq
 \mf d_m(\mf{sl}_n)=\mf d_m(\mf{gl}_n).
 \eeq
 \end{prop}

\begin{proof}
 Following \cite[p.~265, Plate I]{bourbaki4}, we identify $\Lambda'$ with the kernel of the surjection $\Lambda \to \Z$ given by
\beq
(x_1,\ldots, x_n) \mapsto x_1+ \cdots + x_n.
\eeq
Write $\ul{\bf 1}=(1,\ldots, 1) \in \Lambda$. Pick a basis $\mc B'$ of $\Lambda'$; then $\mc B=\mc B' \cup \{ \ul{\bf 1}\}$ is a basis of $\Lambda$. For all $\la \in \Lambda$ we have
\begin{equation} \label{plus 1}
\ms D(\la+ \ul{\bf 1})=\ms D(\la).
\end{equation}
Since $\ms D'$ is the restriction of $\ms D$ to $\Lambda'$, an $m$-period $\varpi$ of $\ms D$ is also an $m$-period of $\ms D'$. Let $\mc F=\mc F_{\mc B}$ and $\mc F'=\mc F_{\mc B}'$. Clearly
\beq
|\mc F|=\varpi \cdot |\mc F'|.
\eeq
It is easy to see using \eqref{plus 1} that
\beq
\mc F  \cap U_m(\ms D)= \amalg_{j=0}^{\varpi-1} \left[j \ul{\bf 1}+ (\mc F' \cap U_m(\ms D'))\right],
\eeq
and therefore
\beq
\frac{|\mc F \cap U_m(\ms D)|}{|\mc F|}=\frac{|\mc F' \cap U_m(\ms D')|}{|\mc F'|}.
\eeq
The equality then follows from Theorem \ref{mega.density.thm}.

\end{proof}

\subsection{Classical Lie algebras} \label{classical.section}
Now let $\mf g$ be a classical Lie algebra, meaning $\mf g=\mf{so}_n$ and $\mf{sp}_{2n}$. (We follow the convention that $\mf{sp}_{2n}$ is a Lie algebra of $2n \times 2n$ matrices.) Let $\rho \in \Lambda^+$ be the half-sum of the positive roots, and put $\mf f(\la)=\ms D(\la- \rho)$. 
The following can be read off from the formulas in \cite[Section 7.1.2]{goodman}:

\begin{prop}
For $\mf g= \mf{so}_n$ and $\mf{sp}_{2n}$, we have
\beq
\mf f \in \mc V_{\mc B} \cap \Q_{\Z}[\Lambda^\vee].
\eeq
Moreover,
\beq
\deg_\bullet(\mf f)=\begin{cases}
2n-1, &  \mf g=\mf{so}_{2n+1} \text{ or } \mf{sp}_{2n} \\
2n-2, &  \mf g=\mf{so}_{2n} \\
\end{cases}.
\eeq
\end{prop}
 
 Proceeding as in the proof of Theorem \ref{main.thm}, we deduce:
 
 \begin{thm} \label{classical.thm} If $\mf g$ is a Lie algebra of type $B_n$, $C_n$, or $D_n$, then
 \beq
 \mf d_m(\mf g)  \leq \omega(m)\exp \left(- \frac{n}{8m} \right).
 \eeq
 \end{thm}

 \begin{proof} (of Theorem  \ref{intro.simple.thm} in the Introduction)
 Except for finitely many  exceptions, every simple Lie algebra is of type $A_n$-$D_n$. Type $A_n$ follows from Proposition \ref{special.prop} and Theorem \ref{main.thm}, while Theorem \ref{classical.thm} covers the other types.
  \end{proof}
 
\subsection{Semisimple Lie Groups}
Let $G$ be a semisimple  Lie group, with Lie algebra $\mf g$.
When $G$ is simply connected, its representations correspond precisely to representations of $\mf g$.
More generally, let $G_{\simp} \to G$ be the simply connected cover, with kernel $C^G$ in the center of $G_{\simp}$, and $\Lambda$ the weight lattice of $\mf g$. 
The irreducible representations of $G$ correspond to the dominant members of $\Lambda^G$, where $\Lambda^G \subseteq \Lambda$ is the sublattice of weights vanishing on $C^G$. 
We have $\Lambda^{G_{\simp}}=\Lambda$, since each  representations of $\mf g$ is the differential of a unique representation of $G_{\simp}$. (See for instance \cite[Section 8.1]{fulton}.) Moreover $[\Lambda: \Lambda^G] =|C^G|$. (In fact, $\Lambda/\Lambda^G$ is in duality with $C^G$ under the perfect pairing described on [Ibid, p.~ 373].)

Put
\beq
\mf d_m(G)=\mf d(\Lambda_m \mid  \Lambda^G);
\eeq
we interpret this as the proportion of irreducible representations of $G$ with degree not divisible by $m$. Of course, $\mf d_m(G_{\simp})=\mf d_m(\mf g)$.
Applying Proposition \ref{sublatticey} gives
\beq
\mf d_m(G) \leq |C^G| \cdot \mf d_m(\mf g).
\eeq
 
Thus an upper bound for $\mf d_m(\mf g)$ leads to an upper bound for $\mf d_m(G)$. For example,
\beq
\mf d_m(\SO_n)  \leq 2 \omega(m)\exp \left(- \frac{n}{8m} \right)
\eeq
and 
\beq
\begin{split}
\mf d_m(\PGL_n) &\leq n \mf d_m(\mf{sl}_n) \\
 &\leq n\omega(m) \exp \left(- \frac{n}{4m} \right).\\
 \end{split}
\eeq

\section{Self-dual Representations} \label{SDR.sec}
Now let $\mf g$ be semisimple. Let $w_0$ be the longest Weyl group element.
Let $\rho^\vee$ be the half-sum of the positive coroots. Put 
\beq
\Lambda_{\sd}=\{ \la \in \Lambda : w_0 \la=-\la \}
\eeq
 and 
\beq
\Lambda_{\orth}=\{ \la \in \Lambda_{\sd} : \lip \la,2\rho^\vee \rip \text{ is even} \},
\eeq
and use the superscript `$+$' to denote the dominant members of these sets.
  According to \cite[Chapter VIII, Section 7, Proposition 12]{bourbaki7}, $\Lambda_{\sd}^+$ is the set of highest weights of irreducible self-dual representations of $\mf g$, and $\Lambda_{\orth}^+$ is the set of highest weights of irreducible orthogonal representations of $\mf g$.    Put
\beq
\mf d_m^{\sd}(\mf g)=\mf d(\Lambda_m \cap \Lambda_{\sd} \mid  \Lambda_{\sd}) \text{ and } \mf d_m^{\orth}(\mf g)=\mf d(\Lambda_m \cap  \Lambda_{\orth} \mid \Lambda_{\orth});
\eeq
we interpret these 
as the proportion of irreducible self-dual, resp., orthogonal representations of $\mf g$ with degree not divisible by $m$.
Since $[\Lambda_{\sd}: \Lambda_{\orth}] \leq 2$, Proposition \ref{sublatticey} gives the bound
\begin{equation} \label{orth.bound}
\mf d_m^{\orth}(\mf g) \leq 2\mf d_m^{\sd}(\mf g).
\end{equation}

For the rest of this section we use material from \cite[p.~265-273, Plates I-IV]{bourbaki4}.

\subsection{Self-dual Representations of $\mf{sl}_n$}
We continue with  the notation of Section \ref{slla.subsec}. The longest Weyl group element $w_0$ acts on $\Lambda'$ by
\beq
(x_1,\ldots,x_n) \mapsto (x_n, \ldots, x_1).
\eeq
The map  $\lambda \mapsto  \la^{\#}$ defined by
\beq
(x_1,\ldots, x_k) \mapsto 
\begin{cases}
&(x_1,x_2, \ldots, x_k,-x_k,-x_{k-1},\ldots, -x_2,-x_1) \text{ if $n=2k$} \\
&(x_1,x_2, \ldots, x_k, 0, -x_k, -x_{k-1},\ldots,  -x_2,-x_1) \text{ if $n=2k+1$} \\
\end{cases}
\eeq
gives an isomorphism of $\Z^k$ with the  sublattice of self-dual weights in $\Lambda'$. 
As usual, let $\rho$ be the half-sum of the positive roots.
Then
\beq
\tilde {\mf f}(\la) =\ms D(\la^{\#}-\rho)
\eeq
 lies in $(\mc V_k) \cap \Q_\Z[x_1, \dots, x_k]$.
Furthermore let $m \in \N$ and $\varpi$ be an $m$-period for $\ms D$. Then for any $\mathbf x \in \Z^k$, 
\beq
\begin{split}
\tilde {\mf f}(\lambda + \varpi \mathbf x) &= \ms D(\lambda^\# + \varpi \mathbf x^\#-\rho) \\
&\equiv \ms D(\lambda^\#-\rho) \mod m \\
&= \tilde {\mf f}(\lambda).
\end{split}
\eeq

So, any $m$-period for $\ms D$ is an $m$-period for $\tilde {\mf f}$. 

\begin{thm} \label{self.dual.bound} For any positive integer $m$, we have 
\beq
\mf d_m^{\sd}(\mf{sl}_n) \leq \omega(m)\exp \left(- \frac{n}{36m} \right)
\eeq
and
\beq
\mf d_m^{\orth}(\mf{sl}_n) \leq 2\omega(m)\exp \left(- \frac{n}{36m} \right).
\eeq
 \end{thm} 

\begin{proof}
Let $r_{q_i}(\ms D)$ be the  periods from Proposition \ref{r.formula}.  As in the proof of Theorem \ref{main.thm}, we have
\beq
    \bm{\varpi}  = \max_i r_{q_i}(\ms D) \leq mn.
\eeq
Now Proposition \ref{cmi.fan}  gives
\beq
\begin{split}
\mf d_m^{\sd}(\mf{sl}_n) &\leq \omega(m)\exp\left(-\frac{k^2}{4 \bm{\varpi} }\right) \\
&\leq \omega(m)\exp\left(-\frac{k^2}{4mn}\right).
\end{split}
\eeq
Since $n \leq 3k$, 
 \beq
\mf d_m^{\sd}(\mf{sl}_n) \leq \omega(m) \exp\left(- \frac{n}{36m} \right).
\eeq
 The second inequality of the theorem then follows from \eqref{orth.bound}.
\end{proof}

 \subsection{Self-dual Representations of Classical Lie Algebras}
 
 All representations of Lie algebras of types $B_n$, $C_n$, and $D_n$ are self-dual, except for $D_n$ when $n$ is odd.
 (This is because the longest Weyl group element $w_0$ acts by $-1$ in these cases.) So, let us consider representations of $\mf{so}_{2n}$ with $n$ odd. 

\begin{thm}
For any odd positive integer $n \geq 3$, we have 
\beq
\mf d_m^{\sd}(\mf{so}_{2n}) \leq \omega(m)\exp \left(- \frac{n}{16m} \right).
\eeq
 \end{thm} 

\begin{proof} In this case, $w_0$ acts on $\Lambda=\Z^n$ by
\beq
(x_1,\ldots,x_n) \mapsto (-x_1, \ldots, -x_{n-1},x_n).
\eeq
Hence, the map  $\lambda \mapsto  \la^{\#}$ defined by
\beq
(x_1,\ldots, x_{n-1}) \mapsto (x_1,x_2, \ldots, x_{n-1},0)
\eeq
gives an isomorphism $\Z^{n-1} \overset{\sim}{\to} \Lambda_{\sd}$. As before, define $\tilde {\mf f}: \Z^{n-1} \to \Z$ by the equation
\beq
\tilde{\mf f}(\la)=\ms D(\la^{\#}-\rho).
\eeq
Then $\tilde{\mf f}$ is $m$-periodic with the same period as $\ms D$. The conclusion follows as in the proof of Theorem \ref{self.dual.bound}. 
\end{proof}

\bibliographystyle{alpha}
\bibliography{refs}
 
\end{document}